\documentclass[letterpaper, 10 pt, conference]{ieeeconf}  

\IEEEoverridecommandlockouts                              

\usepackage{graphicx,xcolor}      
\usepackage{amsmath,amssymb,amsopn}
\usepackage{dsfont}
\usepackage{float}
\usepackage[caption = false]{subfig}
\usepackage{xparse}

\newcommand{\R}{\mathbb{R}}

\newcommand{\one}{\ \mathds{1}}

\newtheorem{theo}{Theorem}
\newtheorem{proposition}{Proposition}

\newtheorem{lemma}{Lemma}
\newtheorem{remarkTheo}{Remark}

\newenvironment{remark}[1][]{\begin{remarkTheo}}{\hfill$\square$\end{remarkTheo}}

\title{\LARGE \bf Dynamic Traffic Reconstruction using Probe Vehicles} 

\author{Matthieu Barreau$^{1}$, Anton Selivanov$^{2}$ and Karl Henrik Johansson${}^1$
	\thanks{$^{1}$ Division of Decision and Control Systems, KTH Royal Institute of Technology Stockholm, Sweden (e-mail: barreau,kallej@kth.se).}%
	\thanks{$^{2}$ Department of Automatic Control and Systems Engineering, The University of Sheffield, United Kingdom (e-mail: a.selivanov@sheffield.ac.uk).}%
	\thanks{The research leading to these results is partially funded by the KAUST Office of Sponsored Research under Award No. OSR-2019-CRG8-4033, the Swedish Foundation for Strategic Research and Knut and Alice Wallenberg Foundation. The authors are affiliated with the Wallenberg AI, Autonomous Systems and Software Program (WASP).}
}

\begin{document}

\maketitle
\thispagestyle{empty}
\pagestyle{empty}

\begin{abstract}
	This article deals with the observation problem in traffic flow theory. The model used is the semilinear viscous Burgers equation. Instead of using the traditional fixed sensors to estimate the state of the traffic at given points, the measurements here are obtained from Probe Vehicles (PVs). We propose then a moving dynamic boundary observer whose boundaries are defined by the trajectories of the PVs. The main result of this article is the exponential convergence of the observation error, and, in some cases, its finite-time convergence. Finally, numerical simulations show that it is possible to observe the traffic in the congested, free-flow, and mixed regimes provided that the number of PVs is large enough. 
\end{abstract}

\section{Introduction}

Traffic congestion is the main source of air pollution and time consumption. Various control schemes have been proposed recently \cite{vinitsky2018lagrangian,vcivcic2018traffic} to deal with congestion. Most of them rely on a full-state knowledge. Here, we design an implementable observer that provides the state estimate.

The quality of the observation highly depends on the model. Traffic flow can be seen as an aggregate finite-dimensional system as in \cite{wang2005real,ferrara2018state}. Such models are useful to describe urban areas but become too complex when large scale traffic flows are studied. To resolve the scalability problem, one can use hyperbolic Partial Differential Equations (PDEs) instead \cite{lighthill1955kinematic}. This gives rise to an infinite-dimensional system \cite{seo2017traffic}. 

For infinite-dimensional hyperbolic systems, boundary observers are usually designed using the backstepping \cite{meurer2013extended}. But, most of the time, it is only possible to observe a specific regime \cite{YU2019183}.
For instance, this prevents the observation of mixed free/congested data which is the most important phase for traffic control. To overcome this issue, using dynamical boundary observers for semilinear hyperbolic systems as in \cite{castillo2013boundary}, \cite{qi2018boundary} proposed a multi-mode observer with an estimation of each linearized mode. The final estimate is a merger of these two observations, leading to a poor state reconstruction.

The previous method deals with a linearized system. This is because it is difficult to draw boundary observers with fixed sensors. Indeed, the characteristic velocity can be positive and negative and, for causal purposes, we need to define which boundary of the road is the input \cite{WRIGHT1973293}. That may lead to a switch between the first and the second boundary condition depending on if the information goes forward or backward to the traffic flow. 

Another kind of boundary observers that corrects this causal problem is with a moving domain. These observers are then with a time-dependent domain delimited by the trajectories of the PVs \cite{WRIGHT1973293,herrera2010incorporation}. This is possible thanks to measurements obtained by mobile phones \cite{amin2008mobile} for instance or well-equipped vehicles. The state reconstruction in this case has been the topic of the recent article \cite{delle2019traffic}. There, the authors provide an algorithm to reconstruct in finite-time the density. This is not a dynamic observation since it uses the characteristic methodology. It is therefore highly subject to noise and disturbance.

The goal of this paper is to derive a dynamical boundary observer using PVs in a multi-mode context, i. e. in congested, free-flow, and mixed regimes. The contributions are:
\begin{enumerate}
    \item the convergence of the observer is proved using Lyapunov arguments, leading to algebraic conditions expressed with matrix inequalities;
    \item a new time-dependent Lyapunov functional, inspired from hyperbolic systems, is proposed for studying the viscous Burgers equation;
    \item the limit case of the inviscid Burgers equation is discussed.
\end{enumerate}
Thanks to these novelties, we aim at reducing the effect of uncertainties in the model or noise in the measurements. 

The organization of the paper is as follows. The second section explains some fundamentals of traffic flow and state the model used throughout the paper. Section~III is devoted to the problem statement. Section~IV studies the observation error and provides a numerically tractable test ensuring the exponential or finite-time convergence.
Then, some numerical simulations are proposed to show the efficiency of the designed observer. The conclusion draws perspectives for future work.

\textbf{Notation:} In this paper, we denote by $\partial_x$ the differentiation along the $x$ axis. A shorter notation is $\partial_x f = f_x$. The functional spaces $H^i(I,D)$ refer to Sobolev spaces from $I$ to $D$ and $C^i$ are the spaces of continuously differentiable functions up to a degree $i$. Similarly, $L^{\infty}(I, D)$ refers to bounded functions. $\R$ stands for the real axis while $\R_+$ is the set of positive real numbers. In the case of symmetric matrices, $\star$ represents the lower diagonal part of the matrix and $A \prec (\preceq) B$ means that $A - B$ is negative (semi-)definite. 

\section{Traffic Flow Model}



\begin{figure}
	\centering
	\includegraphics[width=.47\textwidth]{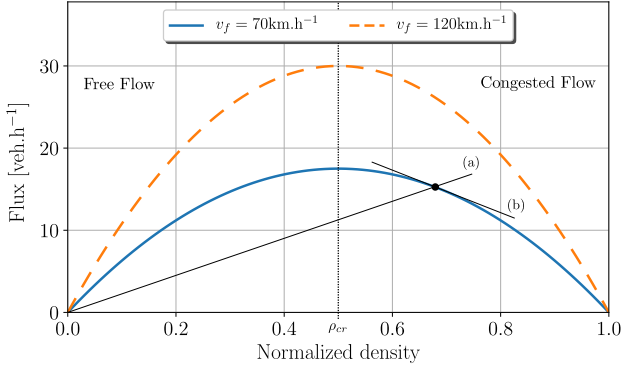}  
	\caption{Fundamental diagram. The slope of (a) is the instateneous speed of the car while the slope of (b) (tangent to the curve) is the characteristic speed.}
	\label{fig:fundDiag}
	\vspace*{-0.7cm}
\end{figure}

Consider the Lighthill--Whitham--Richards model of the traffic flow~\cite{lighthill1955kinematic,richards1956shock}:
\begin{equation} \label{eq:transport}
    \left\{
	\begin{array}{l}
		\partial_t \rho + \partial_x Q(\rho) = 0, \\
		\rho(0,x) = \rho^0(x), 
	\end{array}
	\right.
\end{equation}
where $\rho: \mathbb{R}_+ \times \R \to [0, 1]$ is the density of cars and $Q = \rho V_a$ is the flow with the average speed $V_a(\rho)$. A typical choice is 
\begin{equation} \label{eq:fundamentalDiagram}
    V_a(\rho) = v_f (1 - \rho),\quad Q(\rho) = v_f \rho (1 - \rho), 
\end{equation}
where $v_f$ is the free flow speed. Then the maximum flow corresponds to the critical density $\rho_{cr} = 0.5$. The case $\rho < \rho_{cr}$ corresponds to the free flow, while $\rho > \rho_{cr}$ refers to the congested regime. That leads to the fundamental diagram depicted in Figure~\ref{fig:fundDiag}.

Plugging \eqref{eq:fundamentalDiagram} into \eqref{eq:transport} yields
\begin{equation} \label{eq:burger}
	\partial_t \rho(t,x) + v_f \left( 1-2\rho(t,x) \right) \partial_x \rho(t,x) = 0.
\end{equation}
This is the so-called Burgers equation \cite{burgers2013nonlinear} which is a semi-linear hyperbolic equation and it can be solved explicitly using the characteristic method \cite{evans1998partial}. 
This method is based on the knowledge of iso-density curves $(t, x_c(t))$ such that $\rho(t, x_c(t)) = \rho^0(x_c(0))$. In this case, the PDE transforms into an ordinary differential equation of the form: $\dot{x}_c(t) = 1 - 2\rho(t, x_c(t))$. $\dot{x}_c$ is the analogous of the group velocity $\frac{\partial q}{\partial \rho}$ (slope of (b) in Figure~\ref{fig:fundDiag}) when it comes to study waves. As long as characteristic curves do not intersect, there exists a unique solution to the system with a relatively regular solution \cite{bastin2016stability}, depending on the regularity of $\rho^0$. After an intersection, one needs to use the Rankine-Hugoniot entropy condition \cite{evans1998partial} to ensure the uniqueness of a solution. This solution may then be discontinuous with a shock wave. This shock models the fact that the congestion is moving in time and space. If the shock appears when $\rho \leq \rho_{cr}$, then $\dot{x}_c \geq 0$ and the perturbation is propagating forward in space. Otherwise, $\dot{x}_c \leq 0$ and that means the perturbation goes upstream the traffic flow and backward in space.

The previous paragraph explains that this macroscopic model of traffic flow shows two very different regimes depending on the density on the road. Moreover, it also shows that there might be discontinuities in the solution even if the initial condition is as regular as desired \cite{evans1998partial}.
One way to get around this problem is to add a diffusion term $\gamma > 0$ \cite{whitham2011linear,li2010traffic,nelson2002traveling} as follows:
\begin{equation} \label{eq:viscidBurgers}
	\left\{
	    \begin{array}{l}
	        \partial_t \rho_{\gamma} + v_f \left( 1-2\rho_{\gamma} \right) \partial_x \rho_{\gamma} = \gamma \partial_{xx} \rho_{\gamma}, \\
	        \rho_{\gamma}(0,\cdot) = \rho^0.
	   \end{array}
	\right.
\end{equation}
This new equation is usually referred to as the viscid Burgers equation \cite{burgers2013nonlinear}. Using the Cole-Hopf transformation \cite{hopf1950partial,cole1951quasi}, it is also possible to explicitly solve this equation leading to the following proposition:

\begin{proposition}
	There exists a unique classical solution to \eqref{eq:viscidBurgers} for $\rho^0 \in L^{\infty}(\R, \R)$ and $\rho_{\gamma} \in C^{\infty}((0, \infty) \times \R, \R) \cap C([0, \infty), L^{\infty}(\R, \R))$.
\end{proposition}

Besides the regularization property, it has been proved in \cite{evans1998partial,hopf1950partial} that $\rho_{\gamma}$ converges weakly to $\rho$ when $\gamma \to 0$. We will then use the viscid burger equation as a model and see if the conclusions extend to the case when $\gamma$ goes to $0$.

\section{Observer Design and Preliminaries}

For clarity, this section is divided in two. First, the observer is designed and the problem is stated. Then, some properties of the solution, useful in the sequel, are introduced.

\subsection{Observer design}

We assume that $N$ PVs are located on the road at $x_i \in \R$ for $i \in \{1, \dots, N\}$. These PVs are moving as a particles in the traffic flow, i.e., 
\[
	\left\{ \begin{array}{l}
	    \dot{x}_i(t) = V_a(\rho_{\gamma}(t,x_i(t))) = v_f \left(1 - \rho_{\gamma}(t,x_i(t))\right), \\
	    x_i(0) = x_i^0,
	\end{array} \right.
\]
where $x_1^0 < x_2^0 < \dots < x_N^0$. We assume that PVs can measure the local densities $\rho_{\gamma}(\cdot, x_i(\cdot))$. For each $i \in \{1, 2, \dots, N-1\}$, we introduce an observer
\begin{equation}\label{eq:observer}
	\hspace*{-0.15cm}\left\{
		\hspace*{-0.15cm}\begin{array}{ll}
			\partial_t\hat{\rho}_i + v_f \left(1 - 2 \hat{\rho}_i \right) \partial_x \hat{\rho}_i  = \gamma \partial_{xx} \hat{\rho}_{i},\\
			\hspace{3.5cm} t>0,\quad x\in(x_i(t),x_{i+1}(t)), \\
			\hat{\rho}_i(t,x_{i}(t)) = \rho_{\gamma}(t,x_{i}(t)),\\
			\hat{\rho}_i(t, x_{i+1}(t)) = \rho_{\gamma}(t, x_{i+1}(t)),\\
			\hat{\rho}_i(0,\cdot) = \rho^0(x_{i}^0). 
		\end{array}
	\right.
	\hspace*{-0.3cm}
\end{equation}
The state estimation is given by $\hat{\rho}(t,x) = \hat{\rho}_i(t,x)$ for $x \in [x_{i}(t), x_{i+1}(t)]$. The existence of a solution is dealt by using the Faedo-Galerkin method \cite{evans1998partial}. It consists in introducing a finite-dimensional approximation of $\hat{\rho}$. In a manner similar to the proof of Theorem~2 from \cite{benia2016existence}, one can show that the approximated solution is bounded by its initial condition and consequently converges in a weak-sense. That results in the existence of a unique solution to \eqref{eq:observer} denoted $\hat\rho_i\in H^{1}([0, \infty), H^2(x_{i}, x_{i+1}))$.

Subtracting \eqref{eq:observer} from \eqref{eq:viscidBurgers}, we find that the estimation error $\varepsilon = \rho_{\gamma} - \hat{\rho}$ satisfies 
\begin{equation} \label{eq:error}
	\hspace*{-0.3cm}\left\{
		\hspace*{-0.1cm}\begin{array}{l}
			\partial_t \varepsilon = - v_f\left( 1 - 2\rho_{\gamma}\right) \partial_x \varepsilon + 2 v_f \varepsilon\partial_x \hat{\rho} + \gamma \partial_{xx} \varepsilon, \\
			\hfill t>0,\quad x\in(x_i(t),x_{i+1}(t)), \\
			\varepsilon(t,x_{i+1}(t)) = 0 = \varepsilon(t,x_{i}(t)), \\
            \varepsilon(0,x)=\rho^0(x)-\rho^0(x_{i}^0),\quad x\in(x_i^0,x_{i+1}^0). 
		\end{array}
	\right.\hspace*{-0.3cm}
\end{equation}
We aim to derive the exponential stability conditions for \eqref{eq:error} guaranteeing the existence of $\alpha > 0$ and $K \geq 1$ such that 
\begin{equation} \label{eq:expStab}
	\| \varepsilon(t,\cdot) \|_i \leq K \| \varepsilon(0,\cdot) \|_i e^{-\alpha t}\qquad\forall t \geq 0, 
\end{equation}
where 
\[
    \| \varepsilon(t,\cdot) \|_i^2 = \int_{x_{i}(t)}^{x_{i+1}(t)} \varepsilon^2(t,x) dx.
\]

\subsection{Preliminaries}

This subsection establishes some properties of $\hat{\rho}_i$ that are used in the sequel.


\begin{lemma}[Maximum principle] \label{lem:rhoBounded} The solutions of \eqref{eq:viscidBurgers} and \eqref{eq:observer} satisfy 
\begin{equation}
    \begin{aligned}
    \rho_{\gamma}(t,x) &\in [\inf\rho^0,\ \sup\rho^0],&&\forall t\ge0,\: x\in\R,\\
    \hat{\rho}_{i}(t,x) &\in [\inf\rho^0, \ \sup\rho^0],&&\forall t\ge0,\: x\in[x_i(t),x_{i+1}(t)].
    \end{aligned}
\end{equation}
\end{lemma}
The proof is a combination of the standard maximum principle (see, e.g., \cite{protter2012maximum,clark2011analysis}) and the change of variables $v(y,t) = \hat\rho_i(y(x_{i+1}(t) - x_i(t)) + x_i(t), t)$ with $y\in[0,1]$. 
\begin{lemma} \label{lem:dBounded} If $\inf(\rho^0) > 0$, there exists $d_M > 0$ such that
\begin{equation}\label{dBounded}
0 \leq d_i(t) = x_{i+1}(t) - x_{i}(t) \leq d_M,\quad\forall t \geq 0.
\end{equation}
\end{lemma}
The proof of this lemma is technical and therefore reported to the appendix.

\begin{lemma}[Wirtinger inequality, \cite{hardy1952inequalities}] \label{lem:wirtinger}
If $f\in H^1(a,b)$ is such that $f(a)=0=f(b)$, then 
	\[
		\|f\|_{L^2}\le\frac{b-a}{\pi}\|f'\|_{L^2}.
	\]
\end{lemma}


\section{Convergence Conditions}

In this section, we investigate the exponential stability of the error between the real and observed states.

\subsection{Stability analysis of \eqref{eq:error}}

Here is the main theorem of this article.

\begin{theo} \label{theo:main} Let 
\begin{equation} \label{eq:rhoBounded}
    \rho_{\min}=\inf\rho^0(x)>0,\quad\rho_{\max}=\sup\rho^0(x), 
\end{equation}
and $d_M$ be such that \eqref{dBounded} holds. If there exist $\beta, \xi, p_0 \in\R_+$ and $p_1 \in \R$ such that 
	\begin{equation} \label{eq:opt}
		\Psi_{\xi}(\rho_{\max}) \preceq 0 \quad \text{and} \quad \Psi_{\xi}(\rho_{\min}) \preceq 0,
	\end{equation}
	where
	\[
		\begin{array}{l}
			\Psi_{\xi}(\bar{\rho}) \ = \left[ \begin{matrix}
				\Psi_{11}(\bar{\rho}) + 2 \xi \beta & p_1 - 2 {v_f} \bar{\rho} 
				\\ \star & -2 + p_0 \frac{d_M^2 }{ \gamma \pi^2}e^{\xi \gamma^{-1} d_M}
			\end{matrix} \right], \\
			\Psi_{11}(\bar{\rho}) = {v_f} \left(\rho_{\max} - \frac{8}{3}\bar{\rho} - \frac{4}{3}\rho_{\min} \right) \xi + p_1 \xi + \xi^2
			- \gamma p_0,
		\end{array}
	\]
	then \eqref{eq:error} is exponentially stable. 
\end{theo}

\begin{remark}
The assumption that $\rho^0$ is bounded is reasonable since it holds for $\rho_{\max} = 1$ and we assume the road is not empty. Nevertheless, in the case $\rho_{\max} = 1$, the study presented in this article will be conservative. Since this is a nonlinear equation, it is natural to work locally and therefore consider only parts of the road where $\rho_{\max} < 1$.
\end{remark}


\begin{proof} Due to the zero boundary conditions, the stability of \eqref{eq:error} can be studied independently for each $i\in\{1,\dots,N-1\}$. Consider the functional 
\begin{equation} \label{eq:lyap}
	V(t,\varepsilon(t,\cdot)) = \int_{x_{i}(t)}^{x_{i+1}(t)} \varepsilon^2(t,x)\chi(t,x)dx,
\end{equation}
where 
\begin{equation}
    \chi(t,x)=e^{-\lambda (x_{i+1}(t)-x)}. 
\end{equation}
This functional is inspired by \cite{castillo2013boundary,bastin2016stability}, where hyperbolic systems were considered. We however use it to study a parabolic system. Since $\gamma$ is small, the convection part of \eqref{eq:error} dominates the diffusion part and the "characteristics" propagate backward in space. Thus, the Lyapunov functional has a higher weight towards the right boundary. 

Clearly, 
\begin{equation}\label{LyapCoar}
    e^{-\lambda d_M}\|\varepsilon(t,\cdot)\|_i^2\le V(t,\varepsilon(t,\cdot))\le \|\varepsilon(t,\cdot)\|_i^2. 
\end{equation}
Below we show that the conditions of the theorem guarantee $\dot{V}\le-2 \lambda \beta V$. Indeed, we have 
\[
	\begin{aligned}
		\dot{V}=&- 2 \int_{x_{i}}^{x_{i+1}} v_f \left(1 - 2 \rho_{\gamma}\right) \varepsilon_x\varepsilon \chi+ 4 {v_f}\int_{x_{i}}^{x_{i+1}} \varepsilon^2 \hat{\rho}_x\chi\\
		&+ 2 \gamma \int_{x_{i}}^{x_{i+1}} \varepsilon_{xx} \varepsilon \chi- \lambda v_f \left[1 - \rho_{\gamma}(t,x_{i+1}(t)) \right] V. 
	\end{aligned}	
\]
Since $\chi_x = \lambda \chi$ and $\varepsilon(t,x_{i+1}(t)) = \varepsilon(t,x_{i}(t)) = 0$, integration by parts leads to
\begin{equation} \label{eq:IBP}
	\begin{array}{rl}
		\displaystyle 2 \int_{x_{i}}^{x_{i+1}} \varepsilon_x \varepsilon \chi&= - \lambda V \\
	\end{array}
\end{equation}
and
\[
	\int_{x_{i}}^{x_{i+1}}\varepsilon^2 \hat{\rho}_x\chi= 
		- \lambda \int_{x_{i}}^{x_{i+1}} \varepsilon^2 \hat{\rho}\chi- 2 \int_{x_{i}}^{x_{i+1}} \varepsilon \varepsilon_x \hat{\rho} \chi.
\]
Therefore, 
\[
	\begin{aligned}
		\dot{V}=&\,\lambda \int_{x_{i}}^{x_{i+1}} [{v_f} (\rho_{\gamma}(t,x_{i+1}(t)) - 4 \hat{\rho}) + \lambda \gamma] \varepsilon^2 \chi \\
		&\quad+ 4 v_f \int_{x_{i}}^{x_{i+1}} (\rho_{\gamma} - 2 \hat{\rho}) \varepsilon_x \varepsilon \chi  - 2 \gamma \int_{x_{i}}^{x_{i+1}} \varepsilon_{x}^2 \chi.
	\end{aligned}	
\]
Noting that $\left(\rho_{\gamma} - 2 \hat{\rho}\right) \varepsilon = \varepsilon^2 - \hat{\rho} \varepsilon$, we get
\[
	\begin{aligned}
		\dot{V}=\int_{x_{i}}^{x_{i+1}} \Bigl\{ &\lambda [ v_f (\rho_{\gamma}(t,x_{i+1}(t)) - 4 \hat{\rho}) + \lambda \gamma ] \varepsilon^2 \\
		&+ 4 v_f \varepsilon_x\varepsilon^2 - 4 v_f \hat{\rho} \varepsilon_x \varepsilon- 2 \gamma \varepsilon_{x}^2 \Bigr\}\chi.
	\end{aligned}	
\]
Since 
\[
    \int_{x_{i}}^{x_{i+1}} \varepsilon_x\varepsilon^2 \chi = \frac{\lambda}{3} \! \int_{x_{i}}^{x_{i+1}} (\rho_{\gamma} - \hat{\rho}) \varepsilon^2 \chi,
\]
the previous expression can be written as
\[
	\dot{V} = \int_{x_{i}}^{x_{i+1}} \begin{bmatrix} \varepsilon \\ \varepsilon_x \end{bmatrix}^{\top} \Phi \begin{bmatrix} \varepsilon \\ \varepsilon_x \end{bmatrix}\chi,
\]
where
\begin{equation} \label{eq:phi}
    \begin{aligned}
    &\Phi(t,x) = \left[ \begin{matrix} \phi_{\lambda}(t,x) & -2 v_f {\hat{\rho}(t,x)} \\ \star & -2 \gamma \end{matrix} \right],\\
    &\phi_{\lambda}(t,x) = \lambda v_f \left( \rho_{\gamma}(t,x_{i+1}(t)) - \frac{8}{3} \hat{\rho} - \frac{4}{3} \rho_{\gamma} \right) + \lambda^2 \gamma.
    \end{aligned}
\end{equation} 
Lemma~\ref{lem:wirtinger} and \eqref{eq:IBP} imply 
\[
	\int_{x_{i}}^{x_{i+1}} \left[ \begin{matrix} \varepsilon \\ \varepsilon_x \end{matrix} \right]^{\top} \left(p_0 \Phi_0 + p_1 \Phi_1 \right)  \left[ \begin{matrix} \varepsilon \\ \varepsilon_x \end{matrix} \right] \chi \geq 0
\]
with $p_0 \geq 0$, $p_1 \in \R$, 
\[
	\Phi_0 = \left[ \begin{matrix} -1 & 0 \\ 0 & \frac{d_M^2}{\pi^2} e^{\lambda d_M} \end{matrix} \right] \text{\quad and\quad} \Phi_1 = \left[ \begin{matrix} \lambda & 1 \\ 1 & 0 \vphantom{\frac{d_M^2}{\pi^2} e^{\lambda d_M}} \end{matrix} \right].
\]
Consequently, 
\[
	\dot{V} \leq \int_{x_{i}}^{x_{i+1}} \left[ \begin{matrix} \varepsilon \\ \varepsilon_x \end{matrix} \right]^{\top} \left(\Phi + p_0 \Phi_0 + p_1 \Phi_1 \right) \left[ \begin{matrix} \varepsilon \\ \varepsilon_x \end{matrix} \right] \chi.
\]
The negativity condition $\dot{V} + 2 \lambda \beta V \leq 0$ is equivalent to 
\[
    \Phi(x) + p_0 \Phi_0 + p_1 \Phi_1 \prec \left[\begin{matrix} -2 \lambda \beta & 0 \\ 0 & 0\end{matrix}\right],
\]
which is, by the Schur complement lemma, equivalent to 
\[
    \begin{aligned}
    &-2\gamma + p_0 \frac{d_M^2 }{ \pi^2}e^{\lambda d_M} < 0,\\
    &\phi_{\lambda}(t,x) + 2 \lambda \beta + \lambda p_1 - p_0 + \frac{\left(p_1 - 2 v_f {\hat{\rho}} \right)^2}{2\gamma - p_0 \frac{d_M^2 }{ \pi^2}e^{\lambda d_M}} < 0.
    \end{aligned}
\]
The last inequality holds if 
\[
	\theta_{\xi}(\hat{\rho}) = \bar{\phi}_{\xi}(\hat{\rho}) + 2 \xi \beta - \gamma p_0 + \xi p_1 + \frac{\left(p_1 - 2 v_f \hat{\rho} \right)^2}{2 - p_0 \frac{d_M^2 }{ \gamma \pi^2}e^{\frac{\xi d_M}{\gamma}}} < 0
\]
where $\xi = \gamma \lambda$ and $\bar{\phi}_{\xi}(\hat{\rho}) = \xi v_f(\rho_{\max}-\frac83\hat\rho-\frac43\rho_{\min})+\xi^2$. 
Since $\theta_{\xi}$ is convex, $\theta_{\xi}(\hat{\rho}) \leq 0$ for $\hat{\rho} \in [\rho_{\min}, \rho_{\max}]$ if $\theta_{\xi}(\rho_{\min})\leq0$ and $\theta_{\xi}(\rho_{\max})\leq0$, which are equivalent to $\Psi_{\xi}(\rho_{\min}) \preceq 0$ and $\Psi_{\xi}(\rho_{\max}) \preceq 0$ by the Schur complement lemma. Hence, the conditions of the theorem guarantee $\dot{V} + 2 \lambda \beta V \leq 0$. Using the comparison principle and \eqref{LyapCoar}, we obtain \eqref{eq:expStab} with 
 $K = e^{\xi d_M / (2 \gamma)}$ and $\alpha = \lambda \beta$.	
\end{proof}

For a fixed $\xi$, \eqref{eq:opt} are linear matrix inequalities (LMIs), which can be efficiently solved numerically. The minimum of $\bar{\phi}_{\xi}$ is obtained for 
\[
    \xi^*(\bar{\rho}) = -\frac{v_f\left(\rho_{\max} - \frac{8}{3} \bar{\rho} - \frac{4}{3} \rho_{\min} \right)}{2}.
\]
Since the minimum for $\theta_{\xi}(\bar{\rho})$ is obtained for $\xi \leq \xi^*(\bar{\rho})$, we perform a line-search for $\xi\in[0, \min(\xi^*(\rho_{\min}), \xi^*(\rho_{\max})]$ to verify the conditions of Theorem~\ref{theo:main}. 

\begin{remark}
    The previous stability analysis is conservative since $\Phi(t,x)$ is bounded by a constant matrix. $\Phi_0$ and $\Phi_1$ are introduced to relax this bounding.
\end{remark}

\begin{remark}
    Note that the more PVs there are, the larger $\beta$ can be (and consequently, the decay-rate is larger). Indeed, in this case, $d_M$ decreases and $p_0$ can be larger.
\end{remark}


\subsection{Behavior in the inviscid case}

It has been proved in \cite{burgers2013nonlinear,hopf1950partial} that the solution of \eqref{eq:viscidBurgers} weakly convergence to the solution of \eqref{eq:burger} when $\gamma\to 0$. For $\gamma=0$, the observer \eqref{eq:observer} becomes overdetermined and one of the boundary conditions should be removed. Based on the characteristic methodology, we should remove the left boundary condition to keep the system causal
\begin{equation} \label{eq:observer2}
	\left\{
		\begin{aligned}
			&\partial_t \hat{\rho}_i(t,x)= - v_f \left(1 - 2 \hat{\rho}_i(t,x) \right) \partial_x \hat{\rho}_i(t,x) \\
			&\hat{\rho}_i(t,x_{i+1}(t)) = \rho_{\gamma}(t,x_{i+1}(t)).
		\end{aligned}
	\right.
	\hspace*{-0.3cm}
\end{equation}

If \eqref{eq:opt} stays feasible for a given $\xi > 0$ and $p_0 = 0$, then, by choosing $\lambda = \xi \gamma^{-1}$, we get $\dot{V} + 2 \lambda \beta V < 0$. The error is then exponentially stable for any $\gamma > 0$ and we get
\[
	\forall t \geq 0, \quad 
		\begin{array}[t]{rl} 
			\| \varepsilon(t) \| \!\!\!\!& \leq e^{\lambda \frac{d_M}{2}} e^{-\lambda \beta t} \| \rho^0 - \rho^0(x_i^0) \| \\ 
			& \leq \left(e^{\frac{d_M}{2} - \beta t}\right)^{\lambda} \| \rho^0 - \rho^0(x_i^0) \|.
		\end{array}
\]
Note that $\gamma \to 0$ implies $\lambda \to \infty$, leading to two cases:
\begin{enumerate}
	\item $t \leq t^* = \frac{d_M}{2 \beta}$: $\varepsilon$ stays bounded (cf. Lemma~\ref{lem:rhoBounded});
	\item $t > t^*$: then $\left(e^{\frac{d_M}{2} - \beta t}\right)^{\lambda} \xrightarrow[\gamma \to 0]{} 0$ and $\| \varepsilon(t) \|_i = 0$.
\end{enumerate}
That means $\hat{\rho}$ converges weakly and in finite time $t^*$ to $\rho$.

\begin{remark}
Note that the convergence time $t^*$ is decreasing when there are more PVs (i.e. $d_M$ is decreasing).
\end{remark}

\section{Numerical Simulations and Discussion}

This section is devoted to numerical simulations to conclude on the efficiency of the proposed method.

\subsection{Solution of matrix inequalities~\eqref{eq:opt}}

\begin{figure}
	\centering
	\includegraphics[width=.45\textwidth]{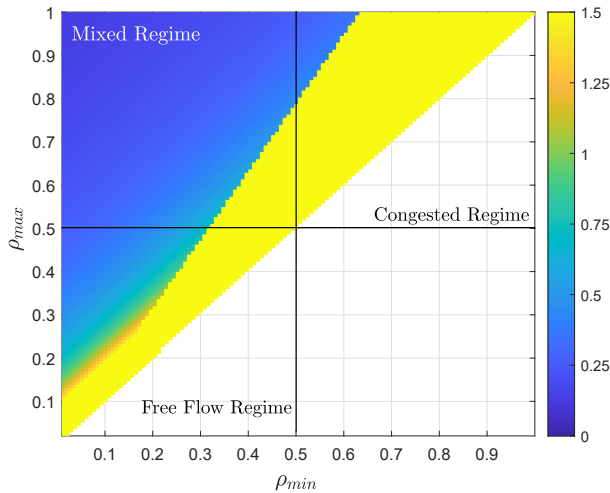}  
	\caption{Maximum $d_M$ for which \eqref{eq:opt} are feasible.}
	\label{fig:rhoMinrhoMax}
	\vspace*{-0.7cm}
\end{figure}

Figure~\ref{fig:rhoMinrhoMax} represents the solution to the optimization problem on a grid of $\rho_{\min}$, $\rho_{\max}$ using a line-search algorithm with $v_f = 70$ km.h${}^{-1}$, and $\gamma = 3$ km${}^2$.h${}^{-1}$. 
The colorbar represents the distance between each PV.

If the maximal distance between two PVs is $0$, then, of course, the observation is possible no matter $\rho_{\min}$ or $\rho_{\max}$. If $d_M > 0$, the quality of the observation decreases as $\rho_{\max}$ increases or $\rho_{\min}$ decreases. That is, the less information there are on the traffic state, the more PVs are needed. 

We can easily see on this plot that the congested regime ($\rho_{\min} \geq 0.5$) is almost everywhere well observed with very few PVs. Since the road is fully congested, the vehicles move at low speed making the observation easier. For the free flow regime ($\rho_{\max} \leq 0.5$), a similar conclusion can be drawn. Indeed, few PVs are enough to observe the road. The main problem in the free flow regime comes from the lack of interaction between vehicles, making the observation sometimes difficult. The mixed regime is the most challenging one and requires to use many PVs.

\subsection{A numerical example}

To show the effectiveness of the main theorem, a numerical simulation is proposed. This is a Godunov-based numerical scheme with $300$ discretization points in space. The proposed example is with $v_f = 70$ km.h${}^{-1}$, $\gamma = 0$ km${}^{2}$.h${}^{-1}$ and has the following initial condition:
\begin{multline} \label{eq:initEx}
	\rho^0(x) = (0.5 + 0.1\sin(5 t)) \one_{(-\infty, 3)} + 0.4 \one_{[3, 3.6)} \\
	+ 0.5 \one_{[3.6, 3.7)} + 0.65 \one_{[3.7, \infty)}.
\end{multline}
The four PVs are located at $x_4^0 = 1.1, x_3^0 = 0.8, x_2^0 = 0.6$ and $x_1^0 = 0.1$. 
\begin{figure}
    \centering
	\subfloat[$\rho$]{\includegraphics[width=.45\textwidth]{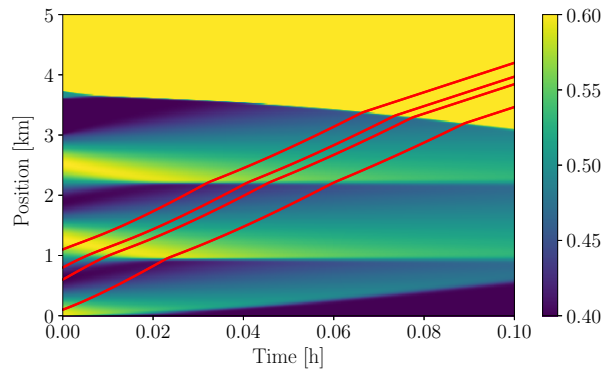}\label{fig:rho}} \\
	\subfloat[$\hat{\rho}$]{\includegraphics[width=.45\textwidth]{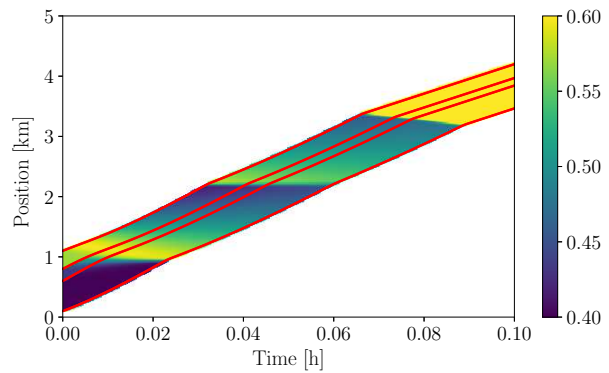}\label{fig:rhoHat}} \\
	\subfloat[Norm of the error between $\rho$ and $\hat{\rho}$.]{\includegraphics[width=.45\textwidth]{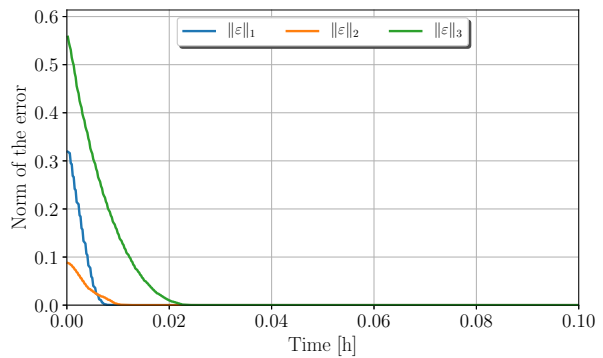}\label{fig:error}}\\
	\subfloat[Distance between PVs.]{\includegraphics[width=.45\textwidth]{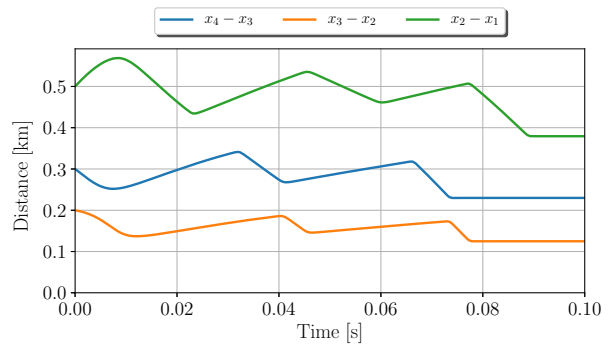}\label{fig:simulationDistance}}
	\caption{Numerical simulation of the original and observed systems.}
	\label{fig:simulation}
	\vspace*{-0.7cm}
\end{figure}
The simulation is in Figure~\ref{fig:simulation} and the trajectories of the PVs are plotted in red. One can see in Figure~\ref{fig:simulationDistance} that the maximal distance between two PVs is $0.6$km and in Figure~\ref{fig:error} that after $1.5$min, the observer has converged. 
However, the matrix inequalities~\eqref{eq:opt} gives a maximal $\beta$ of $0.859$ so a convergence in finite time in less than $20$min. The obtained stability conditions are therefore conservative and we can assess stability of the error only in the case of small $d_M$, even if it is still converging for higher values in simulations.

\section{Conclusion}

In this paper, we derived an observer for traffic flow using probe vehicles. The observer can be seen as a moving boundaries PDE system. Sufficient conditions expressed in terms of linear matrix inequalities with a tuning parameter are given to ensure the exponential or finite-time convergence of the observation error. Numerical simulations show that the convergence area is relatively large. Future work will be devoted to extending this method to other traffic flow models and designing a Luenberger-type observer to increase the stability regions and improve the robustness. 

\appendices

\section{Proof of Lemma~\ref{lem:dBounded}}

    $\bullet$ Note first that $d_i(0) > 0$. Assume that there exists $t$ such that $d_i(t) < 0$. Since $d_i$ is differentiable, that means there exists $t^* \in [0, t]$ such that $d_i(t^*) = 0$ with 
    \[
     \dot{d}_i(t^*) = v_f \left( \rho_{\gamma}(t^*, x_{i}(t^*) - \rho_{\gamma}(t^*, x_i(t^*)^+ \right) < 0.
    \]
    This is impossible since $\rho_{\gamma}$ is continuous. Consequently, we get $d_i \geq 0$.
    
    $\bullet$ Denote by $N_i$ the number of cars between the probe vehicles $i$ and $i+1$. We get the following:
    \[
        N_i(t) = \int_{x_i(t)}^{x_{i+1}(t)} \rho_{\gamma}(t,s) ds.
    \]
    If $\inf(\rho^0) > 0$, then, according to Lemma~\ref{lem:rhoBounded}, we get that:
    \[
        \inf(\rho^0) d_i(t) \leq N_i(t) \leq \sup(\rho^0) d_i(t).
    \]
    The existence of $d_M$ is then ensured if $N_i$ is bounded.
    
    Let us first differentiate $N_i$ along time, that leads to:
    \[
        \dot{N}(t) = \gamma \left\{ \vphantom{\hat\partial_x}\partial_x \rho_{\gamma}(x_{i+1}(t)) - \partial_x \rho_{\gamma}(x_{i}(t)) \right\}.
    \]
    To derive an explicit expression of $\dot N_i$, we can use the Cole-Hopf transformation \cite{hopf1950partial} together with the change of variables defined in \cite{mojtabi2015one}. We then get a solution to \eqref{eq:viscidBurgers}:
    \[
        \rho_{\gamma}(t,x) = \frac{\gamma}{v_f} \frac{w_x(t,x)}{w(t,x)} + \frac{v_f}{2},
    \]
    where $w$ is the solution of the heat equation with diffusion coefficient $\gamma$, written as:
    \[
        w(t,x) = \frac{1}{\sqrt{2 \pi \gamma t}} \int_{-\infty}^{\infty} \phi( \xi ) \exp \left( - \frac{(x - \xi)^2}{\gamma t} \right) d\xi
    \]
    with $\phi$ chosen to satisfies the initial condition.
    Differentiating $\rho_{\gamma}$ yields:
    \[
        \partial_x \rho_{\gamma}(t,x) = \frac{\gamma}{v_f} \frac{w_{xx}(t,x)}{w(t,x)} - \frac{\gamma}{v_f} \frac{w_{x}(t,x)^2}{w(t,x)^2}. 
    \]
    Since $w$ is bounded (maximum principle) and that the integrals exist, one can show that there exists $K$ (which is time-independent) such that 
    \[
        \dot{N}_i(t) \leq K \frac{1}{t^2}.
    \]
    Integrating the previous expression demonstrates that $N_i$ is bounded and that ends the proof.
    
\begin{remark} Note that the condition $\inf \rho^0 > 0$ is not necessary and can be replaced by the improper integral $\int_{-\infty}^{\infty} \rho^0 > 0$. Indeed; as noted in \cite{hopf1950partial}, in this case, there exists $T \geq 0$ such that $\inf \rho_{\gamma}(T,\cdot) > 0$ and the same proof applies straightforwardly.\end{remark}

\bibliographystyle{IEEEtran}
\bibliography{biblio}

\end{document}